\documentclass{amsart}
\usepackage[utf8]{inputenc}
\usepackage{amsfonts}
\usepackage{amsmath}
\usepackage{amssymb}
\usepackage{mathrsfs} 
\usepackage{enumerate}
\usepackage{xcolor}
\usepackage{rotating}
\usepackage{hyperref}
\usepackage{appendix}
\usepackage{chngcntr}
\usepackage[justification=centering]{caption}
\usepackage{tikz}
\usetikzlibrary{decorations.shapes}

\usepackage{lipsum}
\usetikzlibrary{arrows.meta,automata,quotes}
\theoremstyle{plain}
\newtheorem{thm}{Theorem}

\newtheorem{lem}[thm]{Lemma}

\newtheorem{conj}[thm]{Conjecture}

\theoremstyle{definition}

\newcommand{\ignore}[1]{}
\newcommand{\N}{\ensuremath{\mathcal N}}
\renewcommand{\P}{\ensuremath{\mathcal P}}

\newcommand{\n}{\mathbb N}

\renewcommand{\ge}{\geqslant}
\renewcommand{\le}{\leqslant}
\renewcommand{\geq}{\geqslant}
\renewcommand{\leq}{\leqslant}

\usepackage{algorithm,algpseudocode}
\newcounter{algsubstate}

 \usepackage{bbm,dsfont}
 \usepackage{csquotes}

\title{A brief conversation about subtraction games} 
\author{Urban Larsson}
\address{Urban Larsson, IIT Bombay, India}
\email{larsson@iitb.ac.in}
\author{Indrajit Saha}
\address{Indrajit Saha, Kyushu University, Japan}
\email{indrajit@inf.kyushu-u.ac.jp}
\date{\today}
\begin{document}
\maketitle
\begin{abstract}
In this survey we revisit {\sc finite subtraction}, one-heap subtraction games on finite rulesets. The main purpose is to give a general overview of the development, and specifically to draw attention to Flammenkamp's thesis (1997), where he, contrary to other studies, experimentally observes exponential eventual period length of the outcomes, for a carefully selected subclass of games. In addition, we  contribute an appendix on {\sc finite excluded subtraction} by Suetsugu.

\end{abstract}

\section{Introduction}\label{sec:intro}
This survey concerns two-player alternating play impartial combinatorial games, in the normal play convention. The family of ruleset {\sc finite subtraction} is popularly called ``subtraction games''. It is a family of rulesets, played on the non-negative integers. The current player selects a positive number from a given finite ruleset $S\subset \mathbb{N}= \{1,2, \ldots  \}$, and subtracts it from the given non-negative integer. A player who cannot play, because every subtraction leads to a negative number, loses.

 We compile recent developments,  
 by highlighting Flammenkamp's Ph.D. thesis \cite{Flammenkamp_1997}. Contrary to the classical conjecture by Guy, he observes eventual exponential period lengths of the outcomes.  It is interesting to read the quote from Nowakowski's ``Unsolved problems in combinatorial games'' \cite{nowakowski2019unsolved}:
\begin{quote}
``In general, period lengths can be surprisingly long, and it has been suggested that they could be super-polynomial in terms of the size of the subtraction set. However, Guy conjectures that they are bounded by polynomials of degree at most ${n\choose 2}$ in $s_n$, the largest member of a subtraction set of cardinality $n$."
\end{quote}
Let us state here the observation of Flammenkamp as a conjecture. 

\begin{conj}[\cite{Flammenkamp_1997}]
There exists a sequence of subtraction games on finite rulesets with an exponential eventual period length of the outcomes, with respect to the largest member of the subtraction set. 
\end{conj}

In the Appendix, Suetsugu discusses the related {\sc all-but nim}, also known as {\sc finite excluded subtraction}.

\section{Two and three move rulsets.}
\label{sec:whatsub}
 The Children game ``21" starts on a heap of 21 tokens and the two players alternate in removing one or two tokens. Anyone who plays this game a couple of times figures out that they do not want to start. The second player's winning strategy is to complement the other player's move modulo three. However na\"ive this game may seem, together with the game of nim, it is probably the best introduction for an absolute novice to learn about the concept of perfect play and more. 

For another example, if the starting position is $10$ and the subtraction set is $S=\{2,5\}$, then a possible sequence of play is $10\rightarrow 5\rightarrow 3\rightarrow 1$, and the starting player wins. 

The {\em outcome} of an impartial game is a perfect play loss or win, depending on who starts. 
The outcome is \N, if ``the curre\N t  player wins'', and otherwise the outcome is \P, ``the \P revious player wins''. 

Let us compute the initial outcomes $o(\cdot)$ of the ruleset $S=\{2,5\}$. \vspace{2 mm}

{\small
\begin{center}
\begin{tabular}{|c|c|c|c|c|c|c|c|c|c|c|c|c|c|c|c|c|c|}
\hline
$x$ & 0 & 1 & 2 & 3 & 4 & 5 & 6 & 7 & 8 & 9 & 10 & 11 & 12 &13 &14 & 15 & 16 \\ \hline
$o(x)$ & \P  & \P & \N & \N & \P & \N & \N & \P & \P & \N & \N  & \P & \N & \N &\P & \P & \N  \\ \hline
\end{tabular}
\end{center}}\vspace{2 mm}

In the sample play, obviously, the starting player made a mistake, but so did the second player. Similar to the game of ``21", a winning strategy, whenever there is one, is to complement the opponent's move modulo the sum of the available moves.

It is well known \cite{austin1976impartial} that given a 2-move subtraction set $S=\{a,b\}$,  the period length of the outcomes is $a+b$, or $2a$ in case $2a\mid (a+b)$, and there is no preperiod. However, for a 3-move subtraction set $S=\{a,b,c\}$,  the period length of the outcomes is  
more varied. 
As a warm-up, let us compute the initial outcomes for the games $S_1=\{2,3,5\}$, $S_2=\{2,5,7\}$ and $S_3=\{2,4,7\}$, and identify their respective behaviour.\vspace{2 mm}

{\small
\begin{center}
\begin{tabular}{|c|c|c|c|c|c|c|c|c|c|c|c|c|c|c|c|c|c|}
\hline
$x$ & 0 & 1 & 2 & 3 & 4 & 5 & 6 & 7 & 8 & 9 & 10 & 11 & 12 &13 &14 & 15 & 16 \\ \hline
$o_1(x)$ & \P  & \P & \N & \N & \N & \N & \N & \P & \P & \N & \N  & \N & \N & \N &\P & \P & \N  \\ \hline
\end{tabular}
\end{center}}\vspace{2 mm}
The outcomes for ruleset $S_1$ are purely periodic, and the period length is $7=2+5$. 

\vspace{2 mm} 
\begin{center}
\begin{tabular}{|l|l|l|l|l|l|l|l|l|l|l|l|l|l|l|l|}
\hline
$x$    & 0  & 1  & 2  & 3  & 4  & 5  & 6  & 7  & 8  & 9  & 10 & 11 & 12 & 13 & 14 \\ \hline
$o_2(x)$ & $\P$  & $\P$  & $\N$  & $\N$  & $\P$  & $\N$  & $\N$ & $\N$  & $\N$  & $\N$  & $\P$  & $\N$  & $\N$  & $\P$  & $\P$  \\ \hline
$x$    & 15 & 16 & 17 & 18 & 19 & 20 & 21 & 22 & 23 & 24 & 25 & 26 & 27 & 28 & 29 \\ \hline
$o_2(x)$ & $\N$  & $\N$  & $\N$  & $\N$  & $\N$  & $\N$  & $\N$  & $\P$  & $\P$  & $\N$  & $\N$  & $\P$  & $\N$  & $\N$  & $\N$  \\ \hline
\end{tabular}
\end{center}
\vspace{2mm}
The outcomes for ruleset $S_2$ are purely periodic, and the period length is 22. There are many rulesets that start with a preperiod, and, after that, become periodic. One example is $S_3$. \vspace{2mm}

{\small
\begin{center}
\begin{tabular}{|c|c|c|c|c|c|c|c|c|c|c|c|c|c|c|c|c|c|}
\hline
$x$ & 0 & 1 & 2 & 3 & 4 & 5 & 6 & 7 & 8 & 9 & 10 & 11 & 12 &13 &14 & 15 & 16 \\ \hline
$o_3(x)$ & \P  & \P & \N & \N & \N & \N & \P & \N & \N & \P & \N  & \N & \P & \N &\N & \P & \N  \\ \hline
\end{tabular}
\end{center}}\vspace{4 mm}

\begin{figure}[hbp]
 \begin{center}
 \includegraphics[width=\textwidth]{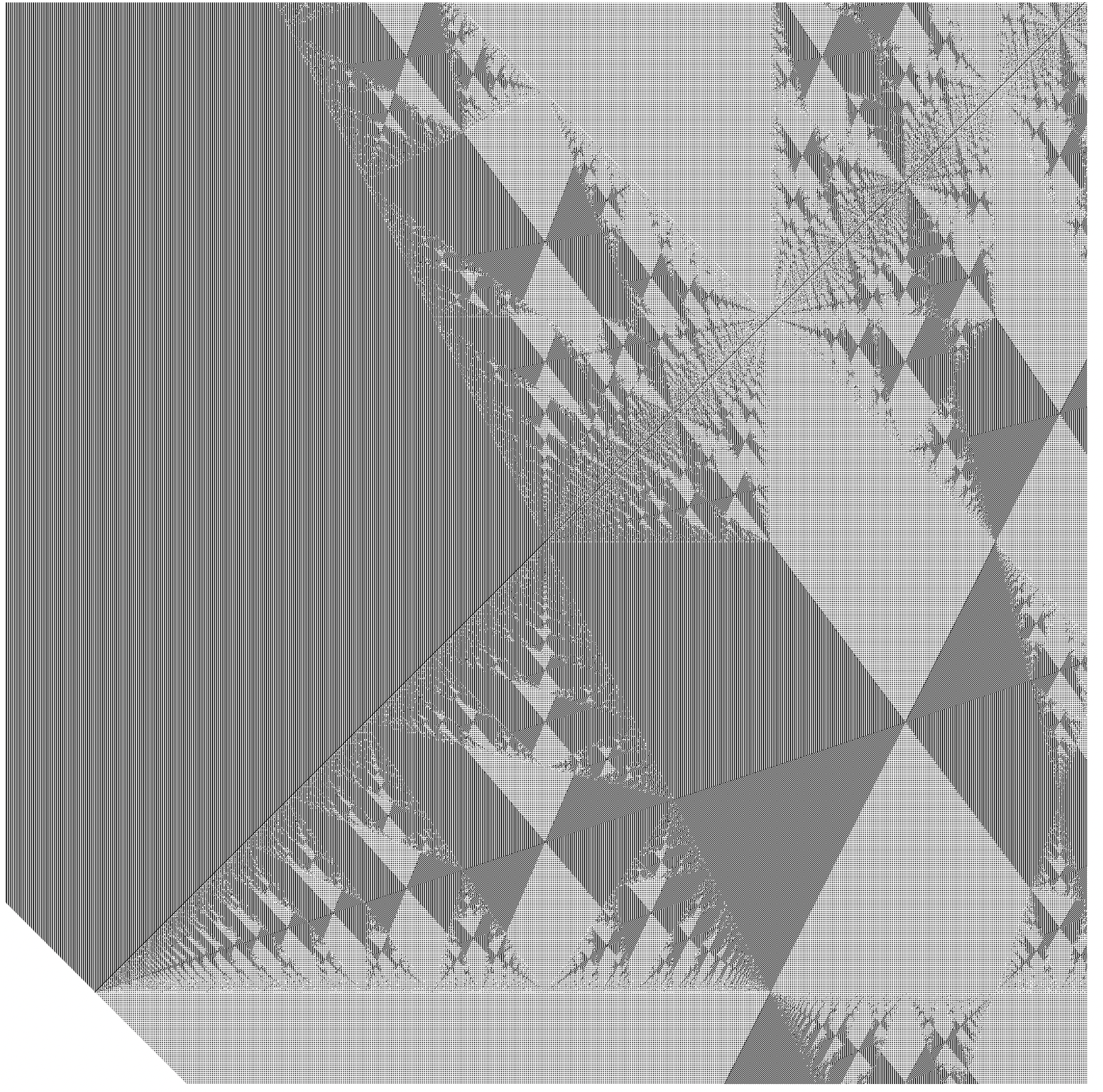}
 \end{center}
 \caption{The picture is adapted from \cite[page 3]{Flammenkamp_1997}, showing a period length classification for 3-move subtraction games, where $s_3=3001, 400<s_1<1600, 900<s_2<2100$. Light gray is period length $p= s_2+s_3$ gray is $p= s_1 + s_3$ dark gray is $p = s_1 + s_2$  
black is the diagonal $s_1+ s_2 = s_3$ and  all other period lengths are white.}
 \label{fig:Flam3move}
 \end{figure}

Why is it that $S_2$ behaves diametric to $S_3$? The latter has a preperiod but a very short period length, whereas the former has no preperiod but a relatively long period length. Is this just a coincidence, and are these types of phenomena chaotic and/or random? Flammenkamp \cite{Flammenkamp_1997} points the direction towards a fractal-like behavior for the classification of 3-move rulesets depending on their period lengths.  He writes: 
\begin{quote}
    ``Die fraktale Struktur, die sich andeutet, scheint in ihrer Form unabh\"angig von $s_3$ und mit wachsendem $s_3$ immer st\"arker ausgepr\"agt zu sein.''
\end{quote} 
That is, the fractal structures increase in detail with increasing $s_3$, and it appears independent of $s_3$. 
For fixed $\max S=s_3$ the fractal-looking behaviors, depends on whether period lengths are of the form $s_1+s_3$, $s_2+s_3$, $s_1+s_2$, or ``something else''. 
See Figure~\ref{fig:Flam3move} for an example.

He shows that very few (a set of density 0) have other period lengths, and that about half of the rulesets have period length $s_1+s_3$. Namely, for increasing $s_3$, the proportions appear to approach $0.53\ldots $, $0.365\ldots $ and $0.105\ldots $  respectively. He also investigates the preperiod lengths of 3-move rulesets and concludes that six forms on a single parameter have the longest pre-periods. Curiously, four out of these satisfy $s_1+s_2  =s_3 \pm 1$. (While, as we will see, 3-move rulesets with the longest periods satisfy $s_1+s_2  =s_3$.) The longest (asymptotic) preperiod appears in the family of games of the form $\{5n-2,5n+3,10n+2\}$, namely, its length is $45n^2-1$, where $n \in \mathbb{N}$.

The recent preprint \cite{mikls2023superpolynomial} has made great progress on 3-move rulesets. We will review this along with the other material.

\section{Some history} 
\label{sec:briefhistory}
A survey on the historical development of combinatorial games appears in \cite{Carvalho2010origin}. 

 According to them, the origin of subtraction games is due to the Italian mathematician  Luca Bartolomeo de Pacioli (written between 1496 and 1508), who proposed an additive version. The game rule says that in an alternating play, a player can add positive integer numbers up to $6$ on their turn, starting with an empty pile, and the goal is to achieve $30$. He provided the solution for the game. Also, in the same survey paper, they mentioned a French mathematician,  Claude Gaspard Bachet de M{\`e}ziriac, who published a set of puzzles in $1612$, and the 22nd problem is similar to  Pacioli's problem, where instead of $30$, a player has to make $100$, and the set of legal moves is to add are any positive integer less than $11$. 
 
 In the 1960s, Golomb  \cite{golomb1966mathematical} defined a generalization of subtraction games to several heaps. 
 His nonlinear shift register model determines the outcomes of any given one-heap subtraction game $S$, with $|S|<\infty$. His construction provides an elementary proof of the ultimate periodicity of any finite subtraction game; it gives an upper bound of the pre-period and period lengths of the form $2^{\max S}$. Indeed, $2^{\max S}$ is the number of combinations of outcome symbols within the range of the options, so a repetition of the outcome patterns must happen.  This upper bound seems very hard to approach, and the big apparent gap with experimental results has triggered a lot of research over the years. This survey aims to shed more light on what has happened since then. 

The nim-value of an impartial game can be computed recursively via a minimal exclusive algorithm on the nim-values of its options (for example, if the option set has nim values $0, 1, 3, 5, 8$, then the position has nim-value $2$). Any \P-position, such as the terminal positions, has nim-value $0$. Nim-sequences are a refinement of outcome sequences because they reveal the winning strategies in disjunctive sum play of impartial normal-play games. In terms of period length, it is easy to see that Golomb's argument instead gives an upper bound of the period length of the nim values of the form $(|S|+1)^{\max S}$. 

In the 1970s, Austin \cite{austin1976impartial} analyses nim-values of subtraction games in his Master's thesis. He demonstrates that a subtraction game has no preperiod if there exists a $p$ such that $p-s\in S$ whenever $s\in S$, and in this case, the period length is (a divisor of ) $p$. He calls such games ``symmetric''. The reason for this is the same as in the children's game ``21'' i.e. $S=\{1,2\}$, where the second player's winning strategy is to complement the first players move modulo $3$. In Austin's generalization, a winning strategy is instead to complement modulo $p$. In more detail, to prove that there is no preperiod, suppose that we have computed all nim values up to a heap of size $p-1$. We must show that the values of $x$ and $p+x$ are the same, say $v$. By induction, we may assume that $y$ and $y+p$ have the same values if $0\le y<x$.  Note that all values smaller than $v$ are available with options to heaps smaller than $x$, and therefore by induction, using the same options from $x+p$, also between $p$ and $p+x-1$. Hence, if no option of $x+p$ has value $v$, then $p+x$ must get this nim-value. Every option is of the form $p+x-s$. Since $S$ is symmetric, this option has an option to $x$. Since the value of $x$ is $v$, for all $s$, therefore $p+x-s\ne v$. 

This idea of complementing or ``reversing out'' the opponent's move (via a known period length) transfers to the following problem: what moves can be adjoined to a given subtraction set without altering its nim-sequence? Austin proves that given any purely periodic subtraction game $S$, with period $p$ say, if $s\in S$ and $p-s>0$, the move $p-s$ can be adjoined to $S$. This is because the second player can reverse the move if desired: if the first player plays $p-s$, the opponent can complement with $s$ to return to the same nim value.   

Without a proof, Austin provides the first complete description of the periods of the nim-sequences of 2-move rulesets, $S=\{a,b\}$, namely, they are  $(0^a1^a)^n2^r$  or  $(0^a1^a)^n0^r2^{a-r}1^r$, if $b=(2n-1)a+r$ or $b=2na+r$, respectively, for some $0<r<a$, $n\in\n$. He claims that the method generalizes to 3-move rulesets in arithmetic progression, i.e. games of the form $S=\{a,b, 2b-a\}$; without a proof, he gives a construction, but it is not very instructive since it lacks game arguments. In support of his construction, he provides the examples $S=\{4,13,22\}$, $\{4,9,14\}$ and $\{4,7,10\}$ with nim-value periods ``$0000 11110000 11112 0000 11112$'', ``$000011110222103321$'' and ``$0000111122223$'' respectively.   We have not yet encountered an update of his result with a proof in the literature. 

Moreover, he studies additive rulesets, that is, rulesets of the form $S=\{a,b, a+b\}$, via what he calls ``the Berlekamp method'' (see also \cite{berlekamp2004winning}). He provides a table of the periods (and preperiods) of all subtraction games with $\max S \le 8$, and he details, for all these games, all moves that can be adjoined without altering the nim-sequences.  Of course, any number for which there is no collision between nim-values can be adjoined as a move. 

Berlekamp et al. \cite[vol. 1, chap. 4, p. 83]{berlekamp2004winning} continues in the early 1980s, and in many instances the study is similar to that of Austin's thesis. For example, they compute the nim-sequences for all rulesets with maximum subtraction $7$ (when Austin already did the same thing up to maximum subtraction $8$). They note that the (ultimate) period length is the sum of two elements from the subtraction set for all but one ruleset, namely the above $S=\{2,5,7\}$. In  \cite[vol. 3, chap. 5, p. 529]{berlekamp2004winning}, they examine additive subtraction games of the form $S = \{a, b, a+b\}$, and provide a complete statement for the nim-sequence whenever $b=2na-r$, for some $0<r<a$, $n\in\n$. Specifically, the period length is $2b+r$. They highlight Ferguson's pairing property to derive the nim-value 1 from the 0s in the harder case when instead $b=2na+r$ for some $0<r<a$, $n\in\n$.\footnote{Ferguson's pairing property states that $\mathcal{G}(n)=1$ if and only if $\mathcal{G}(n-a)=0$, where $a$ represents the smallest number in the subtraction set.} In this case, they claim that the period length is $a(2b+r)$, but, to our best knowledge, the exact nim-sequence has yet been published; see also \cite[Table~4]{ho2015expansion} and the recent preprint \cite{mikls2023superpolynomial}. They omit proofs of these results, and the latter claim still appears open according to \cite{ho2015expansion}; see also recent progress in \cite{mikls2023superpolynomial}.  

Similar to Flammenkamp \cite{Flammenkamp_1997}, Ward \cite{ward2016conjecture} considers 3-move  subtraction $S=\{a,b,c\}$, where $a< b< c$, and conjecture a precise characterization of the nim-value periods. The conjecture is stated in two different regimes. When $S$ is additive, i.e. $c= a+b$, the period is at most a quadratic polynomial in  $a$ and $b$. This part is stated as a theorem in \cite{Flammenkamp_1997}, but concerning only outcomes, and where he cites Winning Ways (where a complete proof appears to be missing); both \cite{austin1976impartial} and \cite{berlekamp2004winning} have at least partial solutions to this part. Both \cite{Flammenkamp_1997} and \cite{ward2016conjecture} conjecture that when $c \neq a+b$, the period length has seven possibilities; in essence, the period length is a divisor of the sum of two elements of the subtraction set.

\section{Polynomial, exponential or something between}
\label{sec:polynomial}
Perhaps the phrase ``it has been suggested'' in the introductory quote is due to  
Alth{\"o}fer  and B{\"u}ltermanns' study \cite{althofer1995superlinear}? They wonder if the parameterized 5-move ruleset $S_a=\{a,8a,30a+1,37a+1,38a+1\}$ has a super-polynomial period length of the outcomes in the parameter $a$.  

The conjecture proposed by Guy suggests that the cardinality of the ruleset is part of the degree of a bounding polynomial. We have not been able to find many examples that are near his bounds for higher cardinalities. However, Flammenkamp \cite{Flammenkamp_1997} conducted extensive experiments, and his experimental observations suggest an exponential period length, even when the cardinality of the subtraction set is a constant, namely five.

Alth{\"o}fer  and B{\"u}ltermann  \cite{althofer1995superlinear} prove some results on the nim-values of the rulesets $S_a= \{a,2a+1, 3a+1\}$, $a \in \n$. This ruleset family has a quadratic polynomial period length in the parameter $a$ and no preperiod. They conjecture that the period length is bounded by a quadratic polynomial in the maximum entry for general 3-move rulesets. 
They prove pure periodicity for $S=\{a,4a,12a+1,16a+1\}$ if $1\le a \le 26$, and where the period length is a cubic polynomial in the parameter $a$, namely $56a^3 + 52a^2 + 9a + 1$. 
At first sight, it appears remarkable that the exact cubic formula holds for all $a\le 26$.  But, in fact, in his thesis,  Flammenkamp claims to have a simpler argument via consistence proof that this result holds for all $a$, and the method extends to other parametrized rulesets. To this purpose, he defines a grammar on words that conjectures periodic nim sequences via for loops, which simplifies the method of proof a lot, leading to fewer sub-cases. Note that Guy conjectured a bounding polynomial of degree 6 for 4-move rulesets. Indeed, Flammenkamp defines a list of 59 parameterized rulesets for which many have period length in terms of a fifth-degree polynomial while some have a sixth-degree polynomial, but none satisfies a seventh-degree polynomial. In this list, he also displays some parametrizations with long pre-periods. Either way, long periods for four move games satisfy one of $s_1+s_2=s_3$, $s_1+s_3=s_4$, $s_2+s_3=s_4$ or $s_1+s_2=s_4$.  To get more wisdom regarding 4-move sets and larger, we continue visiting  Flammenkamp's thesis. Let us summarize the main features as we (amateur translators) understand it; an amateur's translation of selected parts is better than none at all; we have not yet encountered any translation of his amazing work.

Flammenkamp \cite{Flammenkamp_1997} expands on \cite{althofer1995superlinear}. He searches for long outcome periods via extensive computations, and he is not convinced by the mentioned 5-move problem from  \cite{althofer1995superlinear}. He studies instead various rulesets through extensive computations, in some cases via a ruleset parameter or in other cases via a ``record holder'' ruleset in terms of $\max S$.  When it comes to 4-move rulesets he uses the second method to point at some tendency towards exponential vs polynomial behavior. He does this by testing the period lengths of record holders in the range $80\le \max S\le 235$, whether they can be lower bounded by an exponential expression of the form $2^{\alpha\max{S}}$ for some $ \alpha > 0$, or if they can be upper bounded by a monomial of the form $(\max S)^\beta$, for some $\beta>1$. As far as his computation goes, it seems that perhaps both are wrong. Namely $\alpha\sim 0.2$, but with a decreasing tendency, whereas an upper bound $\beta$ seems even more unlikely, because, by indexing with $\max S$,  $\beta_{80}\approx 2.8$ while $\beta_{235}$ has increased to $\approx 5.2$. These record holders for 4-move games all satisfy $s_4-s_3=s_1$. He also lists initial 4-move game record holders with respect to pre-periods. These tend to satisfy $s_1=1$ and $s_4=s_3+1$.  Motivated by these tentative results, Flammenkamp makes more experiments. He completes a table with record holders for any finite ruleset with $\max S\le 30$, and finds a distinguishing property of members of record holder sets. They contain many elements of the form: both $\max S-s$ and $s$ belong to $S$.  Let us define this important property: $S$ is {\em max-symmetric} if, for all $1\le s\le \max S$, if $\max S-s \in S$, then $s\in S$.\footnote{Flammenkamp calls max-symmetric rulesets ``symmetric'' (because he does not study Austin's symmetric games).} By restricting the experiments to max-symmetric sets, he finds that most record holders have size 5. By computing record holders among such sets, for all $61 \le \max S\le 117$, he finds that the period lengths tend to be $2^{\alpha\max S}$, for $\alpha\approx 0.3$. This certainly points towards the existence of subtraction sets with exponential period lengths.

One curious feature, visualized in diagrams in \cite{Flammenkamp_1997}, is that if the subtraction set increases beyond five, then the record holder period lengths tend to rapidly decrease (as a function of $\max S$). Perhaps Guy's conjecture becomes true for large subtraction sets, but remains false for max-symmetric subtraction games of size 5?

Flammenkamp proves that max-symmetric rulesets are purely periodic (Theorem~12, page 42). Note that Austins' argument for symmetric rulesets cannot be used because the move $\max S$ cannot be complemented. Instead, we argue as follows: suppose that we are given an interval of outcomes of length $\max S$, starting at say position $x$. We will compute the outcome of position $x-1$. If, for $s\in S$, one of the positions $x-1+s$ is a \P-position, then clearly, $x-1$ must be an \N-position. Hence, suppose that all positions of the form $x-1+s$ are \N-positions. Since, in particular,  $x-1+\max S$ is an \N-position, it has an option that is a \P-position. But since $S$ is max-symmetric, each option of the form $x-1+\max S-s$, $s\in S\setminus \{\max S\}$, is an \N-position. Hence, the \P-position must be $x-1$. We may assume that our interval comes from the periodic part of the outcome sequence, which leads to a purely periodic pattern. Flammenkamp has already argued for outcomes, but by playing in a disjunctive sum with various nim heaps, we can reuse the argument for nim values and prove that this sequence is also purely periodic. 

Let us list various other results and observations from Flammenkamp's thesis. 
\begin{enumerate}

   \item Conjecture~10 on page~36 is a fairly technical conjecture which describes the periods of 4-move rulesets of the form $\{a,b,2b-a,2b\}$. 
    \item Figure~6 on page~49  illustrates fractal patterns for 5-move max-symmetric rulesets, generalizing those for 3-move games. 
    \item Theorems~15 and 16 on page~50 discuss short period lengths for rulesets that satisfy $s_{i+1}-s_i\le s_1$, for some $i$. If $2s_1\le s_2$ together with another restriction, then the period length is $\le s_3+s_5$. If all moves satisfy the inequality, then the outcome sequence consists of $s_1$ \P-positions followed by $\max S$ \N-positions, etc.  
    \item Observation~17 on page~56: a certain max-symmetric one parameter ruleset with $5n -5$ moves has a longer period than the record holder for max-symmetric ruleset of size 5, if $\max S = 55,70$ or $115$. For example, when $\max S = 115$, then the period length is $604771076188$, compared with $147429129464$ for the record holder, namely $\{15,27,88,100,115 \}$, of max-symmetric rulesets of size 5.
    \item Often, outcomes and nim values have the same period lengths. However, Flammenkamp has found several rulesets that have different period lengths. The following two examples have nim-value period lengths twice that of the outcomes: $\{4,6,11,14\}$ and $\{5,7,14,17 \}$. Further, he describes four 1-parameter rulesets $S_1=\{n, n+3, 3n-1, 3n+3\}$, $S_2=\{2,4n-1,4n+1,4n+5,8n-2\}$, $S_3=\{2,4n+1,4n+3, 4n+7,8n+2\}$ and $S_4=\{n,2n+1,4n+2,5n+3,6n+3\}$ with period lengths $(4n+2,12n+6)$, $(4,8n)$, $(4,8n+4)$ and $(10n+4,10n^2+4n)$ respectively with notation as (outcome period length, nim value period length), pages 57-59, Observations~18 and 20; Corollaries~19 and 21; Theorem~22. He uses his defined grammar to prove these results. (The pre-period lengths are also different.)

    \item Theorem~24 and Corollary~25, page~62. These results significantly lowers the exponential upper bound on period length in case there exist $s_i$ and $s_j$ with $s_i+s_j\le \max S$. Then the period length is upper bounded by $ 2\phi^{\max S}$, where $\phi=\frac{1+\sqrt 5}{2}$. Interestingly, this improved upper bound includes the $\max$-symmetric rulesets for which he experimentally demonstrates exponential period lengths.  
    \item Figure~11 on page~78 illustrates fractal patterns for 4-move rulesets with $s_1+s_3=s_4$, and where $s_3=997$ is fixed.
    \item Flammenkamp's open problems' section: A) Study periodicity properties of rulesets when the initial seed has been altered: any initial configuration of \P~ and \N s of size $\max S$ could initiate computation of periodic `outcomes' (usually this set is only \N~ positions). A game interpretation: if a player moves to a \P ~in this set they win, and otherwise, they lose.  A recent result \cite{mikls2023superpolynomial} shows super-polynomial period lengths even for 3-move rulesets in this generalization.
    B) Study rulesets on the real numbers. In particular, there are infinite rulesets with all moves  $\le \max S$.
\end{enumerate}

Very recently, it has come to our knowledge that the preprint \cite{mikls2023superpolynomial} claims a complete proof of the outcome-periodicity of additive rulesets. We remark that they do not study the nim sequence (which might have a longer period length), as in Austin's thesis, and Winning Ways etc. 

They have made great progress on various outcome sequences of three-move rulesets and in particular: 

\begin{enumerate}
\item  They compute the preperiod and period length of the rulesets $S= \{1,b,c\}$. Ho \cite{ho2015expansion} studied the same ruleset but with the restriction $b< 4c$.

\item As mentioned, Berlekamp et al. \cite{berlekamp2004winning} partially understood the nim sequences of additive three-move rulesets $S=\{a,b,a+b\}$. Mikl{\'o}s and Post\cite{mikls2023superpolynomial} claim a full proof for the outcome sequence. Note that Alth{\"o}fer and B{\"u}ltermann \cite{althofer1995superlinear} solves a specific case of additive games of the form $S= \{a, 2a+1, 3a+1\}$.

\item In traditional normal play subtraction games, one can think of the negative heaps as \N-positions, and where move to a negative terminates the game. Mikl{\'o}s and Post \cite{mikls2023superpolynomial} generalize this idea to a `terminal seed' of any string of outcomes for the positions: $-\max S, \ldots -1$. For example, mis\`ere play would then be the seed $\N^{\min S}\P^{\max S-\min S}$, where the \P s and \N s are interpreted as terminal outcome symbols. In this setting, they claim a proof of superpolynomality of outcome lengths even for three-move rulesets. Namely, they study the parameterized rulesets $S_n=\{n,4n-1,4n^2\}$. By picking the seed starting with a sufficient number of \N -positions, and then the word $\sum_{j=1}^{n-1} \P^j \N^{4n-1-j}$  they claim a proof of superpolynomial period length (but not exponential) in terms of $\max S$, namely their period length is $e^{\Omega(n)}$.

\end{enumerate}

One could imagine doing the same experiment but where the seed is instead a string of nim values. In this way, could one achieve even longer periods?

\section{Related topics}
\label{sec:relatedtopic}
In the previous section, we asked the question: ``which rulesets have long periods, and how long can they be?''. Now, we ask: ``which rulesets have short periods?''. Some rulesets have an ultimate period length of $2$. Obviously, those rulesets cannot have any even number in the subtraction set. 
Cairns and Ho \cite{cairns2010ultimately} define bipartite rulesets: ``an impartial ruleset is bipartite if its game graph is bipartite''. 
They first prove a necessary and sufficient condition for a ruleset to be bipartite. 
     For a finite subtraction set $S$ with $\gcd(S) =1$, the ruleset $S$ is bipartite, for all starting positions, if and only if $1\in S$ and the elements of $S$ are all odd. Of course, the \N-positions are all the odd heap sizes.

They define ultimate bipartite rulesets; a ruleset $S$ is ultimately bipartite if the ultimate period length of outcomes is two, with alternating $\{0,1\}$-nim-values. They provide \cite[Theorem 2]{cairns2010ultimately} three families of ultimately bipartite subtraction sets: 

\begin{enumerate}
    \item $S_k= \{3,5,9, \ldots 2^k +1\}$, for  $k \ge 3$;

    \item  $S_k= \{3,5, 2^k +1\}$, for  $k \ge 3$;

    \item $S_k= \{k,k+2, 2k +3\}$, for  odd $k \ge 3$.
\end{enumerate}

They mention a curious fact about ultimately bipartite rulesets: for large heap sizes, the \N- positions are known, but for small heap sizes winning strategies can be complicated. ``If a ruleset is ultimately bipartite, then,  any sufficiently large heap size $n$ is an \N-position if and only if $n$ is odd.'' The argument is short and elegant. We may assume that $1\not \in S$, for otherwise, we are done.  Suppose, for a contradiction, that a large odd number $n$ is a \P-position. After an even number of moves, the game terminates at an odd number $x<\min S$. Since $x$ is odd, then $x-1\ge 0$, and the same sequence of moves leads to a win for the second player from $n-1$. Hence, the ruleset is not ultimately bipartite, a contradiction.    

Suppose that we adjoin a move to a subtraction set $S$, such that it does not change the nim-sequence. The set of all such elements, including elements in $S$, is called {\em the expansion set} of $S$. 
Roughly, $S$ is non-expandable (according to Ho \cite{ho2015expansion}) if every expansion $s'$ implies that $s'-np\in S$ for some non-negative integer $n$ and where $p$ is the period length. 

Ho \cite{ho2015expansion} continues the work of \cite{berlekamp2004winning} and \cite{austin1976impartial}; he studies the periodicity of 2-move and 3-move rulesets and their expansion sets. He solves the 2-move case completely and displays a number of results for the 3-move case in several tables. For example, he proves the following results on classes of purely periodic games. 
\begin{enumerate}
    \item  Let $a<b$ be relatively prime positive integers, and if $a=1$, then  $b$ is even. 
    Consider the subtraction game $S=\{a,b\}$. If $a+1 < b \le 2a$, then the expansion set is $\{a+1, \ldots, b-1\}+n(a+b)$, $n\in\n$.
    If $b> 2a$, then $S$ is non-expandable. (It follows that if $a=1$ or $a+1=b$, then $S$ is non-expandable.)
    \item Let $a>1$ be an odd positive integer and let $b$ be an even positive integer. The subtraction game $S=\{1,a,b\}$ is purely periodic with  period $a+b$ and nim-values
    $
(01)^{\frac{b}{2}} (23)^{\frac{a-1}{2}} 2.
    $
    The expansion set is     $
\{\{1,3, \ldots , a\} \cup \{b, b+2, \ldots, b+a-1\}\}   +n(a+b), n\in\n$.
\end{enumerate}

He poses a conjecture that ``The subtraction set of a (strictly)  ultimately bipartite game is non-expandable".

Zhang \cite{zhang2024linearity} adjoins a single move $c$ to a subtraction set $S$. He wishes to determine the conditions under which the ruleset $S \cup \{c\}$ exhibits preperiod and period lengths of the outcomes that are linear in the parameter $c$. He provides a conjecture, which states that if $c$ is sufficiently large, and $c$ modulo $q$ is fixed, where $q$ is a multiple of the period length of $S$, then the preperiod and period are linear in $c$. 
He proves this for several cases:  if $S$ contains $1$, if $S$ contains all odd numbers, if $S=\{a, 2a\}$, if $S=\{1,b\}$, or if $S= \{a, a+1, \ldots  b-1, b\}$.

The subtraction game literature and periodicity-related questions expand to partizan subtraction games (where the players have different subtraction sets) \cite{fraenkel1987partizan}. Duchene et al. \cite{duchene2022partizan} compute the period lengths when each player's subtraction set is of size at most two. They show that computing the preperiod length of the outcome sequence is NP-hard in the case where one of the players' subtraction set is of size one. Their method of reduction proof fixes one of the players subtraction set to $\{1\}$. To the best of our knowledge, their method does not generalize to impartial games. And we believe that if both players have non-trivial rules, the hardness problem remains open. 
Golomb's periodicity argument easily generalizes to the partizan setting. But, it turns out that one of the players often dominates the other; for large heap sizes, they win whether they start or not. In the partizan setting, major questions typically concern domination properties and not period lengths etc. 

\section{Subtraction games in more than one dimension}
\label{sec:submore1d}
\begin{figure}[hbp]
 \begin{center}
 \includegraphics[width=\textwidth]{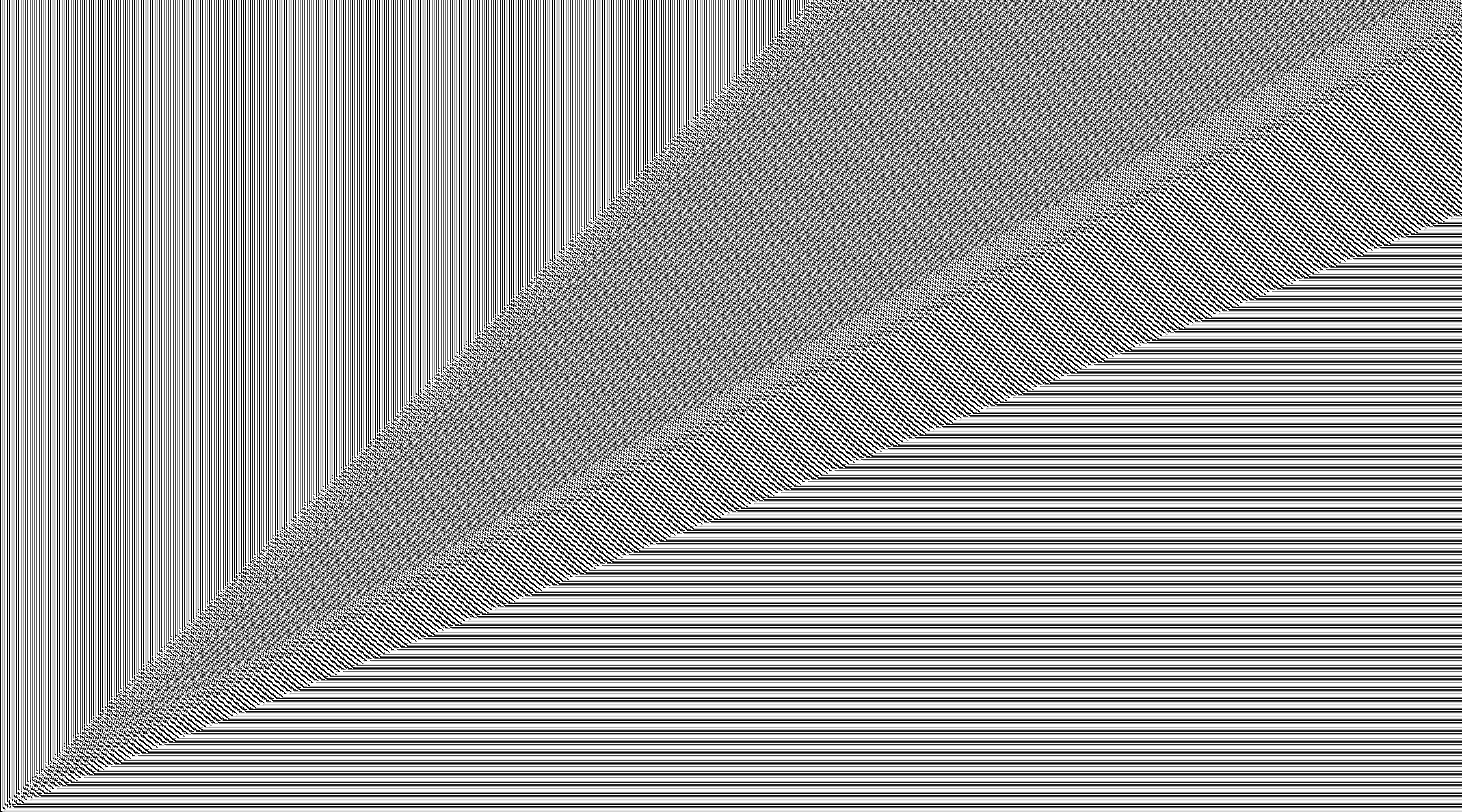}
 \end{center}
 \caption{An outcome visualization of the 2-dimensional  ruleset $S=\{(2,6),(3,3),(6,1),(19,6)\}$ reveals a distinct geometry into five individually regular segments (with bounding strips of constant width). The picture shows the first 3600 by 2000 outcomes, where the \P-positions are colored black. Zoom in for detailed patterns.}
 \label{fig:5segm}
 \end{figure}

It was recently demonstrated \cite{larsson2013heaps} that a slight generalization of subtraction games, studied in more than one dimension, is  Turing complete (see also \cite{fink2012lattice} which proves similar results in a related setting). The dimension and the finite `subtraction set' is taken as input, and the main result states that it is undecidable if two rulesets in the class have the same set of \P-positions. Inspired by this result, the authors focus on play in two dimensions  \cite{larsson2024subtraction}, and by generalizing Golomb's classical method,  prove row and column eventual periodicity of 2-d subtraction games. Moreover, several diagrams expose a frequent  but so far elusive `outcome segmentation' (see Figure \ref{fig:5segm}) 
that lead us to the development of a 0-player game, a ``coloring automaton'' that, in some such cases, can compute the exact set of \P-positions without referring to any \N-position.  We believe that the row-column periodicity result does not generalize to periodicity along any line of rational slope. Several 2-dimensional  rulesets behave irregularly. In particular, the class of (5-move) max-symmetric rulesets from Flammenkamp's thesis generalizes to two dimensions, and they show more complexity than most other rulesets. A conjecture from \cite{larsson2024subtraction} states that, in case of regular outcomes, no 2-dimensional  ruleset on $k$ moves can have more than $k+1$ segments. 

Otherwise the literature on 2-dimensional subtraction games on a finite number of moves is yet remarkably thin. Abuku et al.\cite{abuku2019combination} study a finite version of {\sc cyclic nimhoff}  \cite{fraenkel1991nimhoff}. They provide a closed form for the nim-values whenever the component has a so-called $h$-stair structure.

Carvalho et al. \cite{carvalho2020combinatorics} consider {\sc jenga}, a  popular recreational ruleset. They compute nim-values of {\sc jenga} by viewing them as 2-dimensional  addition games, similar to the ones studied in \cite{larsson2013heaps}. Moreover, they propose a class of impartial rulesets {\sc clock nim} for which {\sc jenga} is a specific example.\\

\noindent {\bf Acknowledgements.} We thank  Koki Suetsugu for contributing the Appendix, Kyle Burke for discussions on NP-hardness and Carlos P. Santos for pointing out the early historical references.  Indrajit Saha is partially supported by JSPS KAKENHI Grant Number JP21H04979 and JST ERATO Grant Number JPMJER2301.

\newpage
\appendix
\numberwithin{thm}{section}
\section{}

\begin{center}
\small\textsc{Koki Suetsugu} \\
\address{Waseda University, Tokyo, Japan}\\
\email{suetsugu.koki@gmail.com} 
\end{center}

The purpose of this appendix is to discuss in detail a variant of the subtraction game {\sc all-but nim} and recent developments on this topic. 
{\sc all-but nim} is studied by Angela Siegel in her Master's thesis \cite{siegel2005finite}. In this variant of  {\sc nim}, the player can remove $n$ stones where $n$ is {\it not} in the given finite set $S \subset \mathbb{N}$. We also refer to the  ruleset as {\sc finite excluded subtraction} and the set $S$ as a FES set.  Siegel proved that the sequence of nim-values of this ruleset is {\em arithmetic periodic}:

\begin{thm}[Siegel \cite{siegel2005finite}]
\label{thm1}
    Let $\mathcal{G}_{\mathbb{N} \setminus S} (n)$ be the nim-value of the position in one-heap {\sc all-but nim} where the FES set is $S$ and the size of the heap is $n$.

    Then, there exist $l \geq 0 , s>0$ and $p > 0$  such that 
    $$
    \mathcal{G}_{\mathbb{N}\setminus S} (n + p) = \mathcal{G}_{\mathbb{N}\setminus S}(n) + s
    $$
    for any $n \geq l$.
\end{thm}
We represent $p$ as the {\em period} and $s$ as the {\em saltus}.
We say the sequence is {\em purely arithmetic periodic} if $
    \mathcal{G}_{\mathbb{N}\setminus S} (n + p) = \mathcal{G}_{\mathbb{N}\setminus S}(n) + s
    $
    holds any $n \geq 0$.

When a sequence is purely arithmetic periodic, it is denoted as $X^n + s$ where $X$ is a sequence and $n$ is the number of repetitions of $X$, and $s$ is the saltus. For example, $(012)^3 +3$ means 
$$
0,1,2,0,1,2,0,1,2,3,4,5,3,4,5,3,4,5,6,7,8,6,7,8,6,7,8,\ldots
$$
Note that the period of the sequence $p$ equals to $|X| \times n$. 

Siegel showed that when $|S| = 1$ or $|S| = 2$, the sequence $\{\mathcal{G}_{\mathbb{N} \setminus S} (n)\}$ is purely arithmetic periodic:

\begin{thm}[Siegel \cite{siegel2005finite}]
If $S = \{a\}$, the sequence of nim-value $\{\mathcal{G}_{\mathbb{N} \setminus S} (n)\}$ is purely arithmetic periodic and the form is 
    $$(0\ldots (a-1))^2  + a.$$

If  $S = \{a, b\} (a<b),$ the sequence of nim-value $\{\mathcal{G}_{\mathbb{N} \setminus S} (n)\}$ is purely arithmetic periodic and the form is 

\begin{eqnarray*}
    \left\{
\begin{array}{ll}
 (0\ldots (a-1))^2  + a & (b \neq 2a)  \\
 (0\ldots (a-1))^3  + a & (b = 2a) \\
\end{array}
\right.
\end{eqnarray*}
\end{thm}

Siegel also studied partizan {\sc all-but subtraction} games in her thesis.

For the case where the size of the FES set is $3$, Siegel conjectured that the sequence of nim-values is purely arithmetic periodic. Later, Sleator and Slusky solve it.

In \cite{sleator2012subtraction}, Sleator and Slusky introduced the FES algorithm, which is a way to calculate nim-values of {\sc all-but nim}. Usually, nim-values are calculated from smaller heaps to larger heaps. In contrast, FES algorithm calculates according to the order of nim-values. That is, the algorithm determines every position which has nim-value $k$ and next it determines every position which has nim-value $k+1$, and so on.

By using this new method, Sleator and Slusky showed another elegant proof of Theorem \ref{thm1}. They also solved Siegel's conjecture:
\begin{thm}[Sleator and Slusky \cite{sleator2012subtraction}]
If  $S = \{a, b, c\} (a<b<c),$ the sequence of nim-values $\{\mathcal{G}_{\mathbb{N} \setminus S} (n)\}$ is purely arithmetic periodic.
\end{thm}

In almost all cases, the sequence of nim-values is equal to the case $S = \{a, b\}, S = \{a, c\},$ or $S = \{a\}$. However, for the case where $S = \{a, b, a+ b\}$, they could not determine the length of the period and introduced a conjecture:

\begin{conj}[Sleator and Slusky \cite{sleator2012subtraction}]
    Let $a$ and $b$ be such that $b > 3a$ and ${\rm gcd}(a, b) = 1$ and let $p$ be the period of the sequence of nim-values $\{\mathcal{G}_{\mathbb{N} \setminus S} (n)\}$, where $S = \{a, b, a+b\}$. If there exists an $m$ that is a multiple of $2a$ with $b < m < a + b$ then 
    $p = 3am$. If no such $m$ exists, then there is some other $n$ with $b < n < a+b$ such that $p = 3an$.
\end{conj}

If $|S| = 4,$ the sequence of nim-values is sometimes non-purely arithmetic periodic. For example, if $S = \{2, 3 ,5, 7\}$, then the sequence of nim-values has a preperiod whose length is $2$.
Abuku and Suetsugu considered the case that the size of FES set is $4$ in \cite{koki2019}.

\begin{thm}[Abuku and Suetsugu \cite{koki2019}]
Let $S = \{a_1, \ldots, a_r\}, g = {\rm gcd}(a_1,\ldots, a_r)$ and $S' = \left \{\dfrac{a_1}{g}, \ldots, \dfrac{a_r}{g} \right\}$. Then, $\mathcal{G}_{\mathbb{N} \setminus S} (n) = \mathcal{G}_{\mathbb{N} \setminus S'} \left(\left \lfloor \dfrac{n}{g} \right \rfloor \right)g + (n \bmod g)$.
\end{thm}
From this theorem, we need to consider only the case that ${\rm gcd}$ of all elements in $S$ is $1$.

\begin{thm}[Abuku and Suetsugu \cite{koki2019}]

If $S = \{a, b, c, d\}(a < b < c < d), $ the sequence of nim-values $\{\mathcal{G}_{\mathbb{N} \setminus S} (n)\}$ is purely arithmetic periodic except for the case $\{a, b, c, d\} = \{a, b, a+b, a+2b\} (b \neq 2a)$ or $\{a, b, c, d\}  = \{a, b, c, a+c\}(b \neq 2a, c \neq 2a)$.

\end{thm}

Also, they calculated the nim-values of the excluded cases and introduced some conjectures.

\begin{conj}[Abuku and Suetsugu \cite{koki2019}]
Consider the case $S= \{a, b, a+b, a+2b\}$. Let $i = \left\lfloor \dfrac{b-a}{2a} \right\rfloor,~ j = \left \lfloor \dfrac{b - (2i + 1)a}{a} \right \rfloor,~ k = b - (2i+j+1)a$ and $ f = 4(i+j)(i+1)a^2 + (4i+3)ka + k^2$. The sequence of nim-values $\{\mathcal{G}_{\mathbb{N} \setminus S} (n)\}$ is purely arithmetic periodic and the period is 

\begin{eqnarray*}
    \left \{
    \begin{array}{ll}
       \dfrac{f}{b-a}  & (a<b<2a)  \\
       \dfrac{f}{\left \lfloor \dfrac{a+b-1}{2a} \right \rfloor a}  & (b \geq 2a, {\rm gcd}(a, b) = a) \\
       \dfrac{f}{{\rm gcd}(a,b)} & (b \geq 2a, {\rm gcd}(a,b)\neq a).
    \end{array}
    \right .\\ 
\end{eqnarray*}

\end{conj}

\begin{conj}[Abuku and Suetsugu \cite{koki2019}]
    The sequence of nim-values $\{\mathcal{G}_{\mathbb{N} \setminus S} (n)\}$ is purely arithmetic periodic if $S = \{a, b, a+b, 2a+b\}$ and $b > 2a$.
\end{conj}

Related to this conjecture, the following lemma has been solved.

\begin{lem}[Abuku and Suetsugu \cite{koki2019}]
    The sequence of nim-values $\{\mathcal{G}_{\mathbb{N} \setminus S} (n)\}$ is purely  arithmetic periodic if $S = \{a, b, a+b, 2a+b\}, b \neq 2a$ and there exists a positive integer $m$  such that $2ma \leq b \leq (2m+1) a$. In this case, the period is $(2m+3)a+b$.
\end{lem}

\begin{conj}[Abuku and Suetsugu \cite{koki2019}]
    Consider the case $S = \{a, b, a+b, 2a+ b\}$. Let $f' = \left \lfloor \dfrac{b+2a-1}{2a}\right \rfloor 4a^2 + 3a(b\bmod a)$. If there exists a nonnegative integer $m$ such that $(2m+1)a <b< (2m+2)a$, then   the sequence of nim-values $\{\mathcal{G}_{\mathbb{N} \setminus S} (n)\}$ is purely arithmetic periodic and the period is $ \dfrac{f'}{{\rm gcd}(a, b)} $.
\end{conj}

In addition, a combination of {\sc all-but nim} and {\sc cyclic nimhoff} is studied in \cite{abuku2019combination}.


\end{document}